\documentclass[11pt]{amsart}
\usepackage{amsmath,amssymb, graphics, amscd,latexsym }

\makeatletter

\newtheorem{Theorem}{Theorem}
\newtheorem{Lemma}[Theorem]{Lemma}
\newtheorem{Corollary}[Theorem]{Corollary}
\newtheorem{Proposition}[Theorem]{Proposition}

\newcommand\modulo{\rm {modulo}\/}
\theoremstyle{remark}
\newtheorem{Remark}[Theorem]{Remark}

\newtheorem{Example}[Theorem]{Example}

\newcommand\Vertex{\rm {Vertex}\/}
\newcommand{\eps}{\varepsilon}

\newcommand\al{\alpha}

\newcommand\si{\sigma}
\newcommand\be{\beta}
\newcommand\Si{\Sigma}
\newcommand\ga{\gamma}
\newcommand\Ga{\Gamma}
\newcommand\de{\delta}
\newcommand\De{\Delta}

\newcommand\cV{\mathcal  V}

\newcommand{\transpose}{{}^{\rm t}\hskip -1pt}
\def\cal#1{{\mathcal #1}}
\def\cC{\mathbb C}

\def\Cone#1{{\rm Cone}\,(#1)}

\def\qed{\hspace*{\fill} $\square$}
\def\QED{\qed}

\begin{document}

\title[Lines on Non-degenerate surfaces 
]
{ Lines on Non-degenerate surfaces  }

\author[
G. Jiang {\tiny and} M. Oka]
{Guangfeng Jiang {\tiny and } Mutsuo Oka}
\address{\vtop{\hbox{Jiang Guangfeng}
\hbox{Department of Mathematics and Physics}
\hbox{Faculty of Science}
\hbox{Beijing University of Chemical Technology}
\hbox{Beijing, 100029}
\hbox{P. R.  China}}
\hspace{2cm}
\vtop{\hbox{Mutsuo Oka}
\hbox{Department of Mathematics}
\hbox{Tokyo Metropolitan University}
\hbox{1-1 Mimami-Ohsawa, Hachioji-shi}
\hbox{Tokyo 192-0397}}}
\thanks{The first author was supported by JSPS: P98028}

\begin{abstract}
On an affine variety  $X$ defined by homogeneous polynomials, every line
in the tangent cone of $X$
is a subvariety of  $X$. However there are many other germs of analytic
varieties  which are not of cone type but  contain ``lines'' passing
through the origin. In this paper, we give a method to determine the
existence and the ``number''  of such lines on
non-degenerate surface singualrities.
\end{abstract}
\maketitle

\pagestyle{headings}


\section{Introduction }\label{S:line}

Let $(X,O)$ be a germ of analytic varieties embedded in $({\mathbb C}^n,O)$
with a singularity at $O$. 
 By abuse of language,  we say that
$L$ is {\it a line} in  $(X,O)$ if 
 $(L,O) $ is
 a smooth curve germ in $(X,O)$ and $L\setminus \{0\}$ is contained
 in the regular part of $X$.

In \cite{J, JS},  lines on hypersurfaces with simple singularities
 are classified
by using the classification machinery. All the hypersurfaces 
of dimension 2 and 3 with
simple or simple elliptic  singularities passing through $x$-axis are
equivalent to (under the coordinate transformation preserving the $x$-axis)
 some surfaces defined by explicit equations. It turns out that the $A,D,E$
singularities  split in this classification. This says that different smooth
curves on the same surface might have different properties.

\vspace{.1cm} 
Let $\pi: \tilde{X} \to (X,O)  $ be a resolution
of a  surface $(X,O)$ with an isolated singularity at the origin $O$ and 
let $\{E_1,\dots, E_r\}$ be the exceptional divisors of $\pi$.
For an exceptional divisor $E_i$, let
${\cal L}_{E_i}$ denote the set of lines  on $(X,0)$ whose  strict transform
intersect  $E_i$ transversally.
It is known that 
$\cal L_{E_i}$ is non-empty if and only if  there exist a function germ $h$
in the maximal ideal 
$\mathfrak m$ such that the multiplicity of $\pi^*h$ along $E_i$ is one
and conversely
 any line in $X$ is contained in some $\cal L_{E_i}$  (\cite{GL1,GL2}).
We call $E_i$  a {\it normally smooth divisor} if 
${\cal L}_{E_i}\ne\emptyset$.
 Geometrically this implies that $d\pi( v)\ne 0$ for any  tangent vector 
$v\in T_{P}\tilde X$ as long as $P\in E_i\setminus\bigcup_{j\ne i}E_j$
and $v$ is not tangent to $E_i$.
If $E_i$ is normally smooth, any germ of a curve
intersecting $E_i\setminus \bigcup_{j\ne i}E_j$ transversely
 defines a line in $X$.
Any two lines in the same 
${\cal L}_{E_i}$ can be connected by an analytic family of lines in $(X,O)$.

\vspace{.15cm}
For a given resolution $\pi:\tilde{X}\to X$,
we consider the integer
  $\rho(\pi):=\sharp \{E_i;{\cal L}_{E_i}\ne \emptyset \}$. 
This number  depends on the  resolution. 
 Put 
${\rho}(X,O)$ to be  the minimal value of $\rho(\pi)$. Obviously 
$\rho(\pi)=\rho(X,O)$
if $\pi:\tilde X\to X$ is a  minimal resolution.
We call $\rho(\pi)$ {\it the line index  of the resolution}
$\pi:\tilde X\to X$
and  we call  $\rho(X,O)$  {\it  the  line index} of $(X,O)$.

\vspace{.15cm}
In this paper, we study $\rho(\pi)$ where $\pi$ is a toric
resolution 
 of a non-degenerate surface
singularity. Let  $(X,0)\subset
({\mathbb C}^3,0)$ be a surface defined by $f(z_1,z_2,z_3)=0$ with
isolated 
 singularity  at the origin. We assume that
$f$ is non-degenerate in the sense of the Newton boundary (\cite{Kou}).
Let  $\Si ^*$ be a regular simplicial cone
subdivision of the dual  Newton diagram $\Gamma ^*(f)$ and let 
$\pi:{X_{\Si^*}}\to (X,0)$
be the associated toric resolution. 
We denote $\rho (\pi) $ by $\rho(\Si^*)$ for simplicity.
To each vertex   $P=\transpose(p_1,p_2,p_3)$ 
of  $\Si ^*$,
there corresponds an exceptional divisor $E(P)$ of $\pi$, which may
have several components. The multiplicity of  $\pi^*z_i$ 
along $E(P)$ is equal to  $p_i$ (\cite{oka}). 
Thus by the result of Gonzalez-Sprinberg
and  Lejeune-Jalabert (\cite{GL1}), $E(P)$ is normally smooth if and only if 
$\min(p_1,p_2,p_3)=1$.
 We  observe  that
$\rho(\Si ^*)$ is independent of the choice of
$\Si ^*$ under certain conditions (see Proposition~\ref{isomorphism}). 
This  allows us to use the  canonical toric resolution to determine
$\rho(\Si^*)$.
Note that
a toric resolution  is not necessarily minimal. So, in general,
 $\rho(\Si ^*)$ may be bigger than $\rho(X,O)$ (see Example~\ref{notminimal}).
However to have the equality $\rho(\Si ^*)=\rho(X,O)$, it is enough that 
$\pi: {X_{\Si^*}}\to X$ is
line-equivalent to
the minimal resolution (see \S~\ref{line-admissible} for the definition).
The purpose of this paper is to give a method to compute
 $\rho(\Si ^*)$.

\section{Line-admissible blowing-ups}\label{line-admissible} 
Let $(X,O)$ be a germ of a surface with an
isolated singularity at $O$.
 Suppose that we have a good resolution $\pi_1:X_1\to
X$ and let $E_1,\dots, E_r$ be the exceptional divisors of $\pi_1$.
Take a divisor $E_{i_0}$ and a 
point $Q$ on $ E_{i_0}$ and let $\pi_Q:\tilde{X}_1\to X_1$
 be the blowing-up
at $Q$ and let $E_Q$ be the exceptional divisor of $\pi_Q$.
The following statements are  obvious.

\begin{Proposition}\label{blowing-up}\sl
Take a function $h\in \mathfrak m$ and let $m_i$ be the multiplicity of 
$\pi_1^*h$ along $E_i$. Then the multiplicity $m_Q$  of the pull-back
$\pi^*_Q(\pi_1^*h)$ along $E_Q$ is the sum of $ m_i$
for all  $i$ such that $Q\in E_i$. In particular,
$m_Q\ge 1$, and $m_Q =1$ if and only if 
$m_{i_0}=1$ and $Q\in E_{i_0}\setminus \bigcup_{i\ne i_0}E_i$.
\end{Proposition}

\begin{Corollary}\label{type2}\sl
Under the situation of Proposition \ref{blowing-up},
$E_Q$ is a normally smooth divisor of the composition 
$\pi_1\circ \pi
_Q:\tilde{X}_1\to X$ if and only if $E_{i_0}$ is a normally
smooth divisor of $\pi_1:X_1\to X$ and $Q$ is  contained in 
$E_{i_0}\setminus \cup_{j\ne i_o}E_j $.
\end{Corollary}

We call $\pi_Q:\tilde{X}_1\to X_1$ a {\it line-admissible}  blowing-up 
if either the center  $Q$ is at the intersection of two exceptional divisor
or the supporting divisor is not normally smooth.
Suppose that we have another good resolution $\pi_2:X_2\to X$.
We say that $\pi_2:X_2\to X$ is {\it line-equivalent} to $\pi_1:X_1\to X$ if 
there exist a finite chain of resolutions
$\pi_i':Y_i\to X,i=1,\dots, s$ such that (1)  $Y_1=X_1$ and $\pi_1'=\pi_1$
and $Y_s=X_2$ and 
$\pi_s'=\pi_2$
 and (2)  any consecutive resolutions factor by 
either $\sigma_i:Y_i\to Y_{i+1}$ or $\si_i':Y_{i+1}\to Y_i$, where 
$\si_i$ and $\si_i'$ are  line-admissible blowing-ups.

An immediate consequence of the definition and Corollary \ref{type2} is:
\begin{Corollary}\label{independent}\sl
Assume that $\pi_i:X_1\to X,i=1,2$  are line-equivalent.
Then $\rho(\pi _1)=\rho(\pi _2)$.
\end{Corollary}

\section{Toric resolution and the computation of $\rho (\Sigma^*)$}

\subsection{Non-degenerate surfaces}\label{non-degeneratesurfaces}
We begin with recalling the toric resolutions of surface singularities
since this also helps us to fix some notations. We use 
the notations of \cite{oka}. Let $(X,O)$ be the germ of 
 a surface in $({\mathbb C}^3,O)$ defined by a function $f:({\mathbb C}^3,O)
\to ({\mathbb C},O)$. 
Hereafter we always assume  that
 $X$ has an isolated singularity at $O$.
Let 
$\sum\limits_{\nu}a_{\nu}z^{\nu}$ be the Taylor expansion of 
 $f$. The {\it Newton polyhedron} $\Gamma _{+}(f)$
is by
definition the convex hull of $\bigcup _{\{\nu; a_{\nu}\ne 0\}}\{\nu +
{\mathbb R}^{3}\}.$ The {\it Newton boundary} 
$\Gamma (f)$ 
is by definition
the union of  the compact faces of $\Gamma _{+}(f)$.

Let $N:={\rm Hom}_{{\mathbb Z}}({\mathbb Z}^{3},{\mathbb Z} )$ be the 
set of covectors. We identify $N$ with ${\mathbb Z}^{3}$ and we denote
the elements of $N$  by 
column vectors. Let 
$N_+$ be the set of covectors
$P=\transpose(p_1,p_2,p_3)\in N$  with $p_i\ge 0, i=1,2,3$.  
 Put $E_1:=\transpose (1,0,0), E_2:=\transpose (0,1,0)$ and 
$ E_3:=\transpose(0,0,1)$.
$P$ is called {\it
strictly positive covector} if $p_j>0$ for all $j$. 
 We denote  
the minimal value of the
linear function $P$
on $\Gamma _{+}(f)$
 by
$d(P; f)$. Put
$\De (P; f)=\{z\in \Gamma_{+}(f)\mid P(z)=d(P ; f)\}$.
 The {\it face function} of $f$ with respect to $P$ is by definition
$f_P(z)=f_{\De (P; f)}
  :=\sum _{\nu \in \De (P; f)}a_{\nu}z^{\nu}$.
Two  covectors $P, P' \in N_+$ are equivalent 
if and only if
$\De (P; f)=\De(P'; f)$.
The {\it  dual Newton diagram} $\Gamma ^*(f) $ of
$X$ is a conical polyhedral subdivision of
$N_+$ given by  the above equivalent classes.

A surface $X$ is called {\it non-degenerate } (with respect 
to the local coordinate $z$) if  for any strictly positive covector
 $P\in N_+$,
$X^*(P):=\{z\in {{\cC}^{*}}^{3}\mid f_P(z)=0\}$
is a reduced
non-singular surface  in the complex torus ${{\cC}^{*}}^{3}.$
The notion of non-degeneracy can be extended to complete intersection
 varieties (cf. \cite{Kho,oka}).

\subsection{Canonical subdivisions}\label{canonicalsubdivision}
We assume that $X$ is defined by $f(z_1,z_2,z_3)=0$  and $f$ is
non-degenerate. 
Let $\Ga^*(f)_2^+$ be the union of the  two-dimensional
cones $\text{Cone}(P,Q)$ of  $\Ga^*(f)$ 
such that the interior points are 
strictly positive. Let $\Sigma^*$ be a regular simplicial subdivision of 
the dual Newton diagram $\Ga^*(f)$ and let
$\pi: X_{\Sigma^*}\to X$ be the associated 
toric modification.
Let $\cV(\Sigma^*)$ be the set of
strictly positive vertices $P$'s of $\Sigma ^*$  such that $\dim\De(P;f)\ge 1$.
The exceptional divisors correspond bijectively to $\cV(\Sigma^*)$ and
for each $P\in \cV(\Sigma^*)$ we denote
the corresponding divisor by $E(P)$. Note that $E(P)$ need not to be
 irreducible
but it is a disjoint union of rational spheres if $\dim \De(P;f)=1$.
The number of connected components is  given by $r(P)+1$, where  $r(P)$
is the number of
integral points on the interior of $\De(P;f)$ (\cite[III\S6]{oka}).
The structure of this resolution $\pi:X_{\Sigma^*}\to X$ depends only on the 
restriction of
$\Sigma^*$ to 
$\Ga^*(f)_2^+$. This follows from the following observation:

\begin{Proposition}\label{isom}\sl
Assume that $\Sigma_1^*$ is a regular subdivision of $\Sigma^*$ such that
$\cV(\Sigma_1^*)=\cV(\Sigma^*)$. Then the canonical morphism 
$\psi: X_{\Sigma_1^*}\to X_{\Sigma^*}$, which is induced by the morphism of the
ambient toric varieties, is an isomorphism.
\end{Proposition}
For any two dimensional cone  $\sigma =\Cone{P,Q}\in \Ga^*(f)$,
there exists a    canonical regular  subdivision of $\sigma$
 which is described as follows.
Denote by $d:=\det (P,Q)$ the
 greatest common divisor of  the absolute values of the $2\times 2$
 minors of the matrix $(P,Q)$. If $d>1$,
there exists a unique integer $d_1, 1\le d_1<d$ such that $Q_1:=(P+d_1Q)/d$
is an integral covector. If $d_1>1$, 
repeat the process for $\Cone{P,Q_1}$, until a regular subdivision of
$\Cone{P,Q}$ is obtained. Let $Q_1,\dots, Q_k$ 
be the covectors obtained  in this way.
 Let 
$d/d_1=[m_1,\dots,m_{\ell}]$ be the continuous fraction expansion.
Then 
$\ell=k$ and the self-intersection number of each component of $E(Q_i)$
is
$-m_i$ (cf. \cite[III]{oka}).
Note that $\De(Q_i;f)=\De(P;f)\cap \De(Q;f)$. This implies $r(Q_i)$ is independent of 
$i=1,\dots,k$ and we denote this number by $r(P,Q)$.
Recall that the continuous fraction is defined inductively
by  $[m_1]=m_1 $ and
$[m_1,m_2,\dots,m_k]=m_1-1/[m_2,\dots, m_k]$.

A regular simplicial cone subdivision of $\Ga^*(f)$ is called a 
{\it  canonical  regular 
subdivision} if its restriction to each cone $\sigma$ in
$\Ga^*(f)_2^+$ is canonical  in the above sense, and we denote it by
 $\Sigma_{\rm can}^*$. The associated toric resolution is called the 
{\it canonical toric resolution } of $X$.

Let $Q=\transpose(q_1,q_2,q_3)$ and $P=\transpose
(p_1,p_2,p_3)$. Put $Q_0=Q$ and $Q_{k+1}=P$
and let   $Q_j:=\transpose (q_{1,j},q_{2,j},q_{3,j}),j=0,\dots,k+1$.
The canonical subdivision enjoys
the following property:

\begin{Lemma}\label{keyLemma}\sl 
Assume that $\Cone{P,Q}\in \Ga^*(f)_2^+$.
Fix an $\ell=1,2,3$.
\begin{itemize}
\item[1)] If $q_{\ell}\le 1$  , then $\{q_{\ell,j}\}_{j=0}^{k+1}$ is
 monotone increasing in
$j$ i.e. $ q_{\ell, j+1}\ge q_{\ell, j}$ for $0\le j\le k$.

\item[2)] If  $q_{\ell}\ge 2$,  then 
either $\{q_{\ell,j}\}$ is monotone increasing or
monotone decreasing in $j$ or  
 there exists a $j_0$
($1\le j_0\le k$) such that  $q_{\ell,j_0}\ge 1$ and 
$$p_\ell=q_{\ell, k+1}\ge \cdots \ge q_{\ell,j_0+1}\ge 
q_{\ell,j_0}
\le q_{\ell,j_0-1}\le \cdots \le q_{\ell,0}=q_{\ell}.$$
\end{itemize}
\end{Lemma}

\begin{proof} 
We prove the assertion 2). If the assertion 
does not hold, there exists an index
$j, 1\le j\le k$ such that 
$q_{\ell,j-1}\le q_{\ell,j}>q_{\ell,j+1}$.
This implies that the self intersection number of each component of 
$E(Q_j)$ is $-(q_{\ell,j-1}+q_{\ell,j+1})/q_{\ell,j}>-2$, which
is a contradiction (cf. \cite[II(2.3) and III(6.3)]{oka}).
The assertion 1) follows from 2) as $Q_j, j=1,\dots, k$ are
strictly positive.
\end{proof}

Let
$\Sigma^*$ be any regular simplicial cone subdivision of $\Gamma^*(f)$ and
 let $\pi :\tilde{X}\to X$ be the corresponding 
toric modification. We denote the line index of  $\pi$ by
$\rho (\Sigma^*)$. Take a two dimensional cone $\sigma=\text{Cone}(P,Q) \in
\Gamma^*(f)_2^+$.
Let
 $Q_0:=Q,Q_1, \ldots , Q_{k }, Q_{k+1}
:=P$ be the canonical subdivision of $\sigma$ 
and let $S_0:=Q, S_1, \ldots , S_{\eta},S_{\eta +1}:=P$ be the vertices of
$\Sigma^*$ on this cone.
 By \cite[II(2.3)]{oka},
$\{Q_0,\dots, Q_{k+1} \}\subset \{S_0,\dots,S_{\eta+1}\}$.
We consider the condition:

\vspace{.15cm}\noindent
($\sharp$):
 $\Sigma^*$ has no vertex in the interior of
$\text{Cone}(Q,Q_1)$.

\vspace{.15cm}
\noindent
We say that $\Sigma^*$ satisfies the ($\sharp$)-condition if it satisfies
($\sharp$)-condition
for any  $\Cone{P,Q}$ in $\Gamma^*(f)_2^+$ such that 
$Q$ is not strictly positive.
The inclusion $\cV(\Sigma_{\rm can}^*)\subset \cV(\Sigma^*)$  implies that 
the following statements.

\begin{Theorem}\label{isomorphism}\sl
There exists a canonical morphism 
$\phi:X_{\Sigma^*}\to X_{\Sigma_{\rm can}^*}$. 
Furthermore  $\phi$  is a composition of line-admissible blowing-ups
 if $\Sigma^*$ satisfies {\rm (}$\sharp${\rm )}-condition.
In particular, the line index $\rho(\Sigma ^*)$ does not depend on the
choice of a toric resolution
associated with any 
regular simplicial subdivision satisfying {\rm (}$\sharp${\rm )}-condition
and $\rho(\Sigma ^*)=\rho(\Sigma _{\rm can}^*)$.
\end{Theorem}
\begin{proof}
Take a two dimensional cone $\sigma=\text{Cone}(P,Q) \in \Gamma^*(f)_2^+$ 
and
  assume that $P$ is strictly positive.
Let
 $Q_0:=Q,Q_1, \ldots , Q_{k }, Q_{k  +1}
:=P$ be the canonical subdivision of $\sigma$ and let
 $S_0:=Q, S_1, \ldots , S_{\eta},S_{\eta +1}:=P$ be the vertices of
$\Sigma^*$
on this cone. Write $S_i=\transpose(s_{1,j},s_{2,j},s_{3,j})$. Assume that
$Q_{i_0}=S_{\nu}$ and
$Q_{i_0+1}=S_\mu$ and $\mu-\nu>1$. Take $S_j$ with $\nu<j<\mu$
and put $\al_j=\det(Q_{i_0},S_j)$ and $\be_j=\det(S_j,Q_{i_0+1})$.
Then $\al_j$ and $\be_j$ are positive integers and 
$S_j=\al_j Q_{i_0+1}+\be_j Q_{i_0}$.
This implies that 
$s_{1,j}>s_{1,\nu}+s_{1,\mu}$. Suppose that 
$s_{1}^{\max}=\max\{s_{1,j}; \nu<j<\mu\}$ and put
$\ga=\min\{\ga;s_{1,\ga}=s_{1}^{\max}\}$.  Then by
\cite[II(2.3)]{oka}
the intersection number of
(each component of)
$E(S_{\ga})$ is $-(s_{1,\ga-1}+s_{1, \ga+1})/s_{1,\ga}>-2$.
Then the negativity of the intersection number implies that 
$s_{1,\ga-1}+s_{1, \ga+1}=s_{1,\ga}$. Thus each component of
$E(S_\ga)$ is  a rational sphere of the first kind.
This implies also  that $S_\ga=S_{\ga-1}+S_{\ga+1}$ and
 $\det(S_{\ga -1},S_{\ga +1})=1$.
Put $\cV'=\cV(\Sigma^*)-\{S_{\ga}\}$.  Then we can extend $\cV'$ to get a
regular simplicial subdivision ${\Sigma^*}'$ such that its restriction to 
$\Ga^*(f)_2^+$ is defined by the vertices $\cV'$.
Thus we get a toric resolution $\pi':X_{{\Sigma^*}'}\to X$.
Changing  $\Sigma^*$ outside of $\Ga^*(f)_2^+$ if necessary, we may assume
by Proposition \ref{isom} that 
$\Sigma^*$ is a subdivision of ${\Sigma^*}'$. Thus we get a canonical
morphism
$\psi:X_{\Sigma^*}\to X_{{\Sigma^*}'}$ which factors $\pi$ by $\pi'$.
By the definition, $\psi$ is the  composition of  blowing-up at
$r(S_\ga)+1$ intersection points of respective 
components of $E(S_{\ga-1})$ and $E(S_{\ga+1})$ in $X_{{\Sigma^*}'}$.
Note that $\psi$ is line-admissible unless $Q$ is not strictly positive and 
 $S_\nu=Q_0$ and $S_\mu=Q_1$. This is the situation where
$\psi$ is the blowing up at the intersection of $E(Q_1)$ and $E(Q)$.
This does not occur if $\Sigma^*$ satisfies $(\sharp)$-condition.
Now the assertion follows by the induction on the cardinality of
$\cV(\Sigma^*)\setminus\cV(\Sigma_{\rm can}^* )$. 
\end{proof}

 \subsection{Computation of $\rho(\Sigma ^*_{\rm can})$}
Let $\pi:X_{\Sigma^*}\to X$ be a toric resolution. 
We assume that $\Sigma^*$ satisfies
the ($\sharp$)-condition. We define
${\mathcal V}_{\rm ns}(\Sigma ^*)
:=\{P\in {\mathcal V(\Sigma^*)}
\mid P \text{ has 1 as  a coordinate  }\}$.
We know that $E(P)$ is a normally smooth divisor if and only if $P\in
\cV_{\rm ns}(\Sigma^*)$.
Thus for each $\text{Cone} (P,Q)\in \Ga^*(f)_2^+$, we define
$\rho_{PQ}:=\#{\mathcal V}_{\rm ns}(\Sigma ^*)
\cap{\Cone{P,Q}}^\circ$, where ${\Cone{P,Q}}^\circ$ is the interior
of  $\Cone{P,Q}$. This number is independent of $\Sigma ^*$
by Theorem~\ref{isomorphism}.
 Recall that $r(P,Q)$ 
is the number of integral points in the interior of $\De(P;f)\cap \De(Q;f)$.
By the definition we have

\begin{eqnarray}\label{formula1}
\rho(\Sigma^*)=\sharp\{P\in \cV _{\rm ns}(\Sigma^*);\dim\De(P;f)=2\}+
\sum_{\Cone{P,Q}\in \Ga^*(f)_2^+}(r(P,Q)+1)\rho_{PQ}
\end{eqnarray}
Thus we need  only to compute $\rho_{PQ}$ for the calculation of
$\rho(\Sigma ^*)$.
Take a cone $\sigma=\Cone{P,Q}$ in
$\Gamma^*(f)^+_2$.
The following gives a practical method to compute $\rho_{PQ}$.
\begin{Theorem}\label{existgeneral}\sl
Let $P=\transpose (p_1,p_2,p_3)$ be strictly positive and let 
$Q=\transpose (q_{1},q_{2}, q_{3})$ and assume that 
$d:=\det(P,Q)>1$. Let $Q_i=\transpose(q_{1,i},q_{2,i},q_{3,i}),i=0,\dots,
k+1$
 be the vertices defining the canonical subdivision from $Q$ with $Q_0=Q$
 and $Q_{k+1}=P$. Fix an $\ell\in\{1,2,3\}$.
Then
\begin{enumerate}
\item For each $1\le i\le k$, there exists positive integers $0<\al_i,\be_i<d$
such that $Q_i=(\be_i P+\al_i Q)/d$.
Putting $\al_0=\be_{k+1}=d,~\al_{k+1}=\be_0=0$, they satisfy the inequality:
\[\al_i>\al_{i+1},\quad \be_i<\be_{i+1},\quad i=0,\dots,k\]
\item  Let $\cV_{\rm ns}^{(\ell)}(P,Q)$ be the set of integral  covectors $R$ 
expressed as 
$R=(\be P+\al Q)/d$ where $\al,\be$ are positive integers satisfying 
\begin{eqnarray}\label{ns-ell}
\begin{cases}&\alpha q_{\ell}+\beta p_\ell=d,~0<\al,\be<d\\
&\alpha q_{k}+\beta p_k\equiv 0\mod d \quad
(k\ne \ell)\end{cases}\label{qis1}
\end{eqnarray}
and let  $\cV_{\rm ns}^{(\ell)}(P,Q;\Sigma_{\rm can}^*)$ be the set of covectors
 $Q_i$, $1\le i\le k$ such that   $q_{\ell,i}=1$. Then 
$\cV_{\rm ns}^{(\ell)}(P,Q)=\cV_{\rm ns}^{(\ell)}(P,Q;\Sigma_{\rm can}^*)$.
Note that the inequality $\al,\be<d$  follows automatically
 from the positivity if 
both $p_\ell$ and $q_\ell$ are  positive.
\end{enumerate}
\end{Theorem}
\begin{proof}
The first assertion  follows by an inductive argument.
Write $Q_i=(\be_i P+\al_i Q)/d$ with positive rational numbers $\al_i,\be_i$.
As $\det(P,Q_i)=\al_i$ and $\det(Q_i,Q)=\be_i$, $\al_i,\be_i$ are
positive integers. By the definition of $Q_1$,
 we can write $Q_1=(P+\al_1 Q)/d$ for some $0<\al_1<d$.
The assertion for $Q_1$ holds and $\det(P,Q_1)=\al_1$.
Assume that 
$Q_j=(\beta_j P+\al_j Q)/d$ with $0<\al_j<d$. As $\det(P,Q_j)=\al_j$
and $\{Q_j,\dots,Q_{k+1}\}$ is the vertices of the canonical subdivision
of $\Cone{P,Q_j}$, there exists $\alpha'$, $0<\alpha'<\al_i$, such that 
\[Q_{j+1}=\frac{1}{\al_j}P+\frac{\al '}{\al_j}Q_j=
\frac{1}{\al_j}P+\frac{\al'}{\al_j}\frac{(\beta_j P+\al_j Q)}d
=(\frac{1}{\al_j}+\frac{\al'\be_j}{\al_jd})P+
\frac{\al'}d Q
\]
Thus $\al_{j+1}=\al'<\al_j<d$. 
The inequality $\be_{j+1}>\be_j$ can be proved similarly by using the fact that
$\{P,Q_k,\dots,Q_1,Q\}$ is the vertices of the canonical subdivision of 
the cone $\Cone{P,Q}$ from $P$ (cf. \cite[II(2.3)]{oka}).
Now we show the second assertion. 
The inclusion 
$\cV_{\rm ns}^{(\ell)}(P,Q;\Sigma_{\rm can}^*)\subset
 \cV_{\rm ns}^{(\ell)}(P,Q)$ is obvious.
Suppose that $R=(\be P+\al Q)/d\in \cV_{\rm ns}^{(\ell)}(P,Q)$ 
is not contained in 
$\cV_{\rm ns}^{(\ell)}(P,Q;\Sigma_{\rm can}^*)$.
 Suppose that $R\in \Cone{Q_i,Q_{i+1}}^\circ$.
Then we can write $R=mQ_i+nQ_{i+1}$ for some positive integers $m,n$.
If $i\ge 1$, this gives a contradiction by comparing the $\ell$-th coefficient:
$1=m q_{\ell,i}+nq_{\ell,i+1}\ge m+n$.
Suppose that $i=0$. Write 
$Q_1=(P+\al_1 Q)/d$ as above. Then 
$R=mQ+(P+\al_1 Q) n/d= n P/d+(md+n\al_1) Q/d$.
Thus we get $\al=md+\al_1 n\ge d$ which contradicts to the assumption.
\end{proof}
\begin{Remark}
The computation of $\cV_{\rm ns}(P,Q)$ is most difficult
 for  the  case $p_\ell,q_\ell>1$. Assume that $p_\ell,q_\ell>0$.
If  we have a solution $(\al_0,\be_0)$,
the other solutions 
are reduce to the following equation.
Put $\al=\al_0+\al', \be=\be_0+\be'$. Then 
\begin{eqnarray}
\begin{cases}&\alpha' q_{\ell}+\beta' p_\ell=0\\
&\alpha' q_{k}+\beta' p_k\equiv 0\mod d \quad
(k\ne \ell)\end{cases}
\end{eqnarray}
Let $\Delta:=\Delta(P;f)\cap\Delta(Q;f)$.
Let $T=\transpose(t_1,t_2,t_3)$ be a covector in $\cV_{\rm ns}^{(\ell)}(P,Q)$
(thus $t_\ell=1$).
Geometrically this implies that  $\Delta(T;f)=\Delta$.
In particular, 
$\Ga_+(f)\subset \{(\nu_1,\nu_2,\nu_3); t_1\nu_1+t_2\nu_2+t_3\nu_3\ge
d(T;f)\}$. 
This gives a practical way to find such a $T$.

The case $q_\ell=0$ or $1$, the computation is much easier. See Corollary
\ref{existspecial}.
\end{Remark}

The canonical subdivision of $\Cone{P,Q}$ takes sometimes 
a lot of computations (see Example~\ref{example1}).
Theorem \ref{existgeneral}  gives us a criterion on the
existence or non-existence of normally smooth divisors, without 
computing the whole subdivision
$Q_i,i=1,\dots, k$.
\begin{Example}\label{example1}
 For simplicity, we write
$x=z_1,y=z_2,z=z_3$.
 Let us consider $f(x,y,z)=x^m+y^n+x^ry^r+z^2$.
We assume that $m,n>2r$.
Put $n=n_1r+n_0, m=m_1r+m_0$ with $0\le m_0,n_0\le r-1$.
Then $\Ga(f)$ has two compact faces whose covectors are
$P=\transpose(2(n-r),2r,nr)/\delta_1$ and $Q=\transpose(2r,2(m-r),mr)/\delta_2$
where $\delta_1=\gcd(2(n-r),2r,nr)$ and $\delta_2=\gcd(2r,2(m-r),mr)$ and the corresponding 
dual Newton diagram is as in Figure~\ref{Fi:1}.
Note that 
$d:=\det(P,Q)$ is given by 
$d=2(mn-mr-nr)/(\delta_1\delta_2)$.
We consider $\cV_{\rm ns}^{(1)}(P,Q)$. 
First we consider the covector $T_0=\transpose(1,1,r)$, which is
a weight vector of 
$x^ry^r+z^2$. As  $m,n>2r$, $T_0$ must be on $\Cone {P,Q}$.
To
proceed the further computation, let us assume that 
$n,m,r$ are odd and  $\gcd(m,r)=\gcd(n,r)=1$. This implies $\delta_1=\delta_2=1$.
By Theorem~\ref{existgeneral}, we have
$$
\begin{cases}
&2\be(n-r)+2\al r=d\\
&2\be r+2\al (m-r)\equiv 0  \mod d\\
&\be nr+\al mr \equiv 0  \mod d
\end{cases}
$$
First  we have a canonical solution
$(\alpha_0,\beta_0)=(n-2r,m-2r)$ which corresponds to
the covector $T_0=\transpose(1,1,r)$.
Thus putting $\alpha=\alpha_0+a$ and $\beta=\beta_0+b$, we can reduce the equation as
$$
\begin{cases}
&2b(n-r)+2a r=0\\
&2b r+2 a (m-r)\equiv 0  \mod d\\
& b nr+ a mr \equiv 0  \mod d
\end{cases}
$$
Taking the positivity of $\alpha, \beta$ into account, we have the 
solution
\[
\{(\alpha,\beta)\}=\left\{\left ((n-2r)+2j
(n-r),(m-2r)-2jr\right);
 0\le j\le \left [\frac{m_1-2}{2}\right]\right\}
\]
For example, consider the easiest case $m=n$.
This has a unique solution $(\al,\be)=(n-2r,n-2r)$ and 
$\cV_{\rm ns}^{(1)}(P,Q)=\{B\}$ where $B=\transpose(1,1,r)$. By symmetry, we have
$\cV_{\rm ns}^{(2)}=\{B\}$. 
Note $r(P,Q)=1$. By writing down the equation described by
Theorem~\ref{existgeneral},  we can show $\cV_{\rm ns}^{(3)}(P,Q)=\emptyset$.
\begin{figure}[h]
\setlength{\unitlength}{1pt}
\begin{center}
\scalebox{0.7}
{\includegraphics{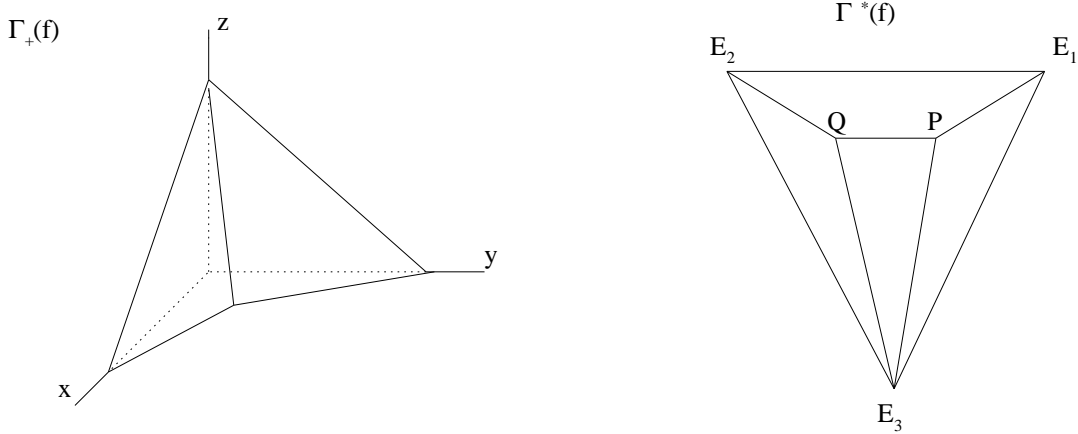}}
\caption{The Newton polyhedron and the dual Newton diagram}
\label{Fi:1}
\end{center}
\end{figure}

Now we look at $\Cone{P,E_1}$ and $\Cone {P,E_3}$. Note that
$\det(P,E_1)=r$ and $\det(P,E_3)=2$.
 It is easy to see that there are no normally smooth 
divisor on these  cones.
Observe that the computation of canonical subdivision of $\Cone{P,Q}$ is 
not so easy.
For example, if $r=15,n=37$,
then  $B=\transpose(1,1,15)$ and first covector $B_1$ (from $Q$) is given by 
$(P+223 Q)/518=\transpose(13,19,240)$ and 
$518/223=[3, 2, 2, 12, 2, 2, 3]$
and it takes some computation to complete the subdivision.
\end{Example}

The following lemma  describes the covectors corresponding to
the  non-compact faces.

\begin{Lemma}\label{non-compact-face}\sl
Assume that $X=\{f(z_1,z_2,z_3)=0\}$ and assume that $f$ is  non-degenerate and  
  $\Ga(f)$ has at least one compact
two dimensional  face for simplicity.
Suppose that $z_2=z_3=0$ is a line in $X$. (So $f$ is not convenient.)
 Then there
is a unique covector
$Q=\transpose(q_1,q_2,q_3)\in \Vertex (\Ga^*(f))$ such that $q_1=0$. Furthermore 
$Q$ takes the form $\transpose(0,1,q_3)$ or $\transpose(0,q_2,1)$.

There exists a unique covector $P=\transpose(p_1,p_2,p_3)$ which corresponds to a
compact divisor and adjacent to $Q$ in $\Ga^*(f)_2^+$.
Then we have $\det(P,Q)=p_1$.
\end{Lemma}
\begin{proof}
As $X$ has an isolated singularity, $f$ must contain a monomial of type
$z_1^a z_2$ or $z_1^az_3$. Suppose that $B:=(a,1,0)\in \Ga(f)$.
Let $C=(b,0,c)$ be  the vertex of $\Ga(f)\cap \{z_2=0\}$ adjacent to
$B$ by an edge.
It is clear that the non-compact face $\Xi$ which has $\overline{BC}$ 
as a face and is unbounded to the direction of the $z_1$-axis
 has covector $Q=\transpose(0,c,1)$. 
One can see that there exists no other 
non-compact face which is unbounded to the  $z_1$-axis direction and bounded to
$z_2,z_3$-direction.
 Let $\De$ be
the compact face which has $\overline{BC}$ as a boundary and let 
$P=\transpose(p_1,p_2,p_3)$  be the corresponding covector.
As $\De(P;f)$ contains $B,C$, we need to have
$ p_1a+p_2=bp_1+cp_3$. 
Now the last assertion follows from 
$\det(P,Q)=\gcd(p_1,p_2-cp_3)=\gcd(p_1,p_1(b-a))=p_1$.
\end{proof}

The following corollary describes explicitly
 $\cV_{\rm ns}^{(1)}(P,Q)$ in the case
$q_1=0$ or $1$.
\begin{Corollary}\label{existspecial}\sl
With the assumptions of Theorem~\ref{existgeneral},   we have the following.
\begin{itemize}
\item[1)] Assume $q_{1}=0$. Then 
$\cV_{\rm ns}^{(1)}(P,Q)\ne \emptyset$  if and only if $d:=\det(P,Q)>1$
and $d=p_1$.
In this cases, $\cV_{\rm ns}^{(1)}(P,Q)=\{Q_1\}$.
 If $Q\ne E_2,E_3$, 
  then  $\{y=z=0\}\subset X$ and   $d=\det(P,Q)=p_1$.

\item[2)] Assume $q_{1}=1$. Then $\cV_{\rm ns}^{(1)}(P,Q)\ne \emptyset$
 if and only if
$d >p_1$. In this case, %
 we have 
$Q_i=(iP+(d-ip_1)Q)/d$ for $i=1,\dots,[d/p_1]$ and 
$\cV_{\rm ns}^{(1)}(P,Q)=\{Q_i;i=1,\dots, [d/p_1]\}$.
\end{itemize}
\end{Corollary}

\begin{proof}
Assume that $Q'=(\be P+\al Q)/d\in \cV_{\rm ns}^{(1)}(P,Q)$ with $0<\al,\be<d$.

1) If $q_1=0$,  we have 
$\gcd(q_2,q_3)=1$. As
$d=\gcd(p_1q_2,p_1q_3,p_2q_3-p_3q_2)=\gcd(p_1,p_2q_3-p_3q_2)$, 
$d$ divides $p_1$. Thus $Q'\in  \cV_{\rm ns}^{(1)}(P,Q)$ if and only if $d=p_1$ and
$\be=1$. In this case, $Q'=Q_1$ and $\cV_{\rm ns}^{(1}(P,Q)=\{Q_1\}$. 
Assume that $Q\ne E_2,E_3$.
By the definition of $\Ga^*(f)_2^+$, $\De(Q;f)$ 
is a non-compact face with dimension 2. In particular, $\{y=z=0\}\subset X$.
By Lemma \ref{non-compact-face},
we have  $d= p_1$. 

2) Suppose that $q_1=1$. Then $\be p_1+\al=d$. This implies
$d>p_1$. Put $d=rp_1+d'$ with $0\le d'<p_1$ and $r=[d/p_1]$. 
Then by the above equality, we have
$(\al,\be)=(d-j p_1,j),~j=1,\dots,[d/p_1]$. 
Put $Q_j':=(jP+(d-jp_1)Q)/d$.
By the definition, $d$ divides the minors of $(P,Q)$ which are 
$p_1q_2-p_2,p_1q_3-p_3,p_2q_3-p_3q_2$. 
Thus 
$\be p_j+\al q_j=\be p_j+(d-\be p_1)q_j\equiv \be(p_j-p_1q_j)\equiv 0
\mod d$ for $j=2,3$. Thus $Q_j'$ is an integral covector
for $\be=1,\dots,r$. 
It is clear that $Q_1'=Q_1$. Assume that  
$Q_r'=Q_\iota$ for some $\iota$.
By the monotonity of the coefficients (Lemma \ref{keyLemma}),
we have $Q_j\in \cV_{\rm ns}^{(1)}(P,Q)$ for $j\le \iota$. 
Thus $\iota=r$ and $Q_j'=Q_j$ for $j\le r$.
\end{proof}
\begin{Remark}
In the case of non-convenient surface with $q_1=0$, 
the divisor $E(Q_1)$
corresponds to the  deformations of the line $z_2=z_3=0$. In fact,
$E(Q)$ is a non-compact divisor which is the strict transform of
$z_1$-axis and $E(Q)$ intersects transversely with $E(Q_1)$.
\end{Remark}
For $R\in \cV_{\rm ns}^{(\ell)}$,  write  $R=(\be P+\al Q)/d$. 
We call $\be/d$ 
{\it the P-coefficient} of $R$.
\begin{Corollary}\label{Theta1}\sl
With the assumptions of Theorem~\ref{existgeneral}, suppose that $q_{1}>1$. 
Let $\bar Q=(\bar \be P+\bar \al Q)/d\in \cV_{\rm ns}^{(\ell)}$
and $\underline{Q}=(\underline {\be} P+\underline{ \al} Q)/d\in \cV_{\rm ns}^{(\ell)}$
be the covectors with maximal and minimal $P$-coefficients in $\cV_{\rm ns}^{(\ell)}$.
Then 
\begin{equation}\label{qis3}
\rho_{PQ}^{(\ell)}=1+|\det(\bar Q,\underline{Q})|
=1+\frac{|\bar\be \underline{\al}-\bar\al \underline{\be}|}{d}
\end{equation}
\end{Corollary}
\begin{proof}
Denote by  $d ':=|\det(\bar Q,\underline{Q} )|$. 
Suppose that $\underline{Q}=Q_i$ and $\bar Q=Q_{i+j}$.
Then $\cV_{\rm ns}^{(\ell)}=\{Q_i, \dots,  Q_{i+j}\}$ by Lemma \ref{keyLemma}
and $\rho_{PQ}^{(\ell)}=j+1$.
By the assumption, we  have $Q_{i+1}=(Q_{i+j}+(d'-1)Q_i)/d'$.
As the continuous fraction $d'/(d'-1)$ is given by  $[2,\dots,2]$
($(d'-1)$  copies of 2), we get $j-1=d'-1$ and the assertion
follows  immediately.
\end{proof}

\section{Applications} \label{application}

\subsection{Weighted homogeneous surfaces}\label{all-types}
In this section we study lines on  weighted
homogeneous surface singularities, which are classified 
as follows ( \cite{OW,oka}):

$\begin{array}{ll}
X_{\rm I}:&h_{\rm I}=x^a+y^b+z^c=0,\cr
X_{\rm II}:&h_{\rm II}=x^ay+y^b+z^c=0,\cr
X_{\rm III}:&h_{\rm III}=x^ay+xy^b+z^c=0,\cr
X_{\rm IV}:&h_{\rm IV}=x^ay+y^bz+z^c=0,\cr
X_{\rm V}:&h_{\rm V}=x^ay+y^bz+z^cx=0,\cr
X_{\rm VI}:&h_{\rm VI}=xy+z^c=0,\cr
X_{\rm VII}:&h_{\rm VII}
=x^{a}z+y^{b}z+z^{c}+tx^{c_1}y^{c_2}=0,\quad t\ne 0\cr
X_{\rm VIII:}&h_{\rm VIII}
=x^{a}y +xy^{b}+xz^{c}+ty^{c_1}z^{c_2}=0, \quad t\ne 0.
\end{array}
$

\vspace{.2cm}
The surface $X_{\rm I}$ is called a Pham-Brieskorn surface.
This type of surfaces have been studied in the previous paper
\cite{JOPS}.
The surface $X_{\rm VI}$ is an $A_{c-1}$ type singularity.
There are exact $c-1$ families of lines on this surface 
(see \cite{GL1,GL2,JOPS, JS}).
On surface $X_{\rm VII}$ and $X_{\rm VIII}$, the term $y^{c_1}z^{c_2}$
must be on the supporting plane of the previous three monomials.
Thus $a,b,c$ are not arbitrary.
The Newton
boundaries of the  surfaces other than $X_{\rm VI}, X_{\rm VII}$
  and $X_{\rm VIII}$ are triangles. 
Note that for a  weighted homogeneous surface, the Newton boundary
has only one compact  2-dimensional face.
Let $P=\transpose(p_1,p_2,p_3)$ be the corresponding covector.
The formula (\ref{formula1}) in \S\ref{line-admissible} reduces to
\begin{eqnarray}\label{formula1'}
\rho(\Sigma^*_{\rm can})=\varepsilon+
\sum_{\Cone{P,Q}\in \Ga^*(f)_2^+}(r(P,Q)+1)\rho_{PQ}(\Sigma^*_{\rm can}).
\end{eqnarray}
where $\varepsilon=1$ if $P\in \cV_{\rm ns}(\Sigma_{\rm can}^*)$  and 
$\varepsilon=0$  otherwise.

For each type of surfaces,  one can calculate
 $\rho_{PQ}(\Sigma^*_{\rm can})$
for each $\Cone{P,Q}$ in the dual Newton diagram by using the method
described in the previous sections.

\begin{Lemma}\label{numberPEi}\sl 
Assume that $\Cone{P,E_i}$ be  a cone in 
$\Ga ^*(f)_2^+$. Then $\det(P,E_i)$ is given by  $\delta
_i:=\gcd(p_j, p_k)$ with
$\{i,j,k\}=\{1,2,3\}$.
Assume that  $\delta _i>1$.
\begin{itemize}
\item[1)] $ \cV_{\rm ns}^{(i)}(P,E_i)\ne \emptyset$
          if and only if  $\delta _i>p_i$ and 
$\rho_{PE_i}^{(i)}=\left[\frac{\delta _i}{p_i}\right]$.
\item[2)] $\cV_{\rm ns}^{(j)}(P,E_i)\ne \emptyset $
 if and only if $p_j|p_k$. In  this case, $\rho_{PE_i}^{(j)}=1$.
\item[3)] 
\[
\rho_{PE_i}=\begin{cases}
&0,\quad\text{if}~ \left[\frac{ \delta
_i}{p_i}\right]=0 ~\text{and}~\de_i<\min\{p_j,p_k\} \\
&\max\{1,\left[\frac{
\delta _i}{p_i}\right]\},\quad\text{otherwise}\end{cases}\]
   \end{itemize}
\end{Lemma}
\begin{proof}
This follows from Corollary~\ref{existspecial}.
\end{proof}

\begin{Lemma}\label{numberPQ}\sl Let $\Cone{P,Q}$ be a cone in 
$\Ga ^*(f)_2^+$ with   $Q=\transpose(0,c,1)$. Suppose that
$\det(P,Q)=p_1>1$. Then 
\[
\rho_{PQ}=
\begin{cases}
&\max\{1,\left[\frac{p_1}{p_2}\right],\left[\frac{p_1}{p_3}\right]\},
\quad c=1\\
&\rho_{PQ}^{(2)}+\max\{1,\left[\frac{p_1}{p_3}\right]\}-\eps,\quad
c>1
\end{cases}
\]
where $\eps=1$ if either $Q_1\in \cV_{\rm ns}^{(2)}(P,Q)$ or 
$Q_{j_1}\in \cV_{\rm ns}^{(2)}(P,Q)$ 
with $j_1:=\left[\frac{p_1}{p_3}\right]\ge 1$
and $\eps=0$ otherwise.
\end{Lemma}
\begin{proof} Let $Q_1, \ldots, Q_{k}$ be the primitive covectors in
$\Cone{P,Q}$ inserted by the canonical subdivision from $Q$.
If $c=1$, the assertion is immediate from Corollary \ref{existspecial},
as
$q_{1,1}=1$.
We assume that $c>1$. If $\left[{p_1}/{p_3}\right]=0$, the assertion is
obvious.
Assume that $\left[{p_1}/{p_3}\right] \ge 1$.
By Corollary \ref{existspecial}, $Q_j$ is given
by $(jP+(p_1-jp_3)Q)/p_1$ for $1\le j\le j_1$. Thus $q_{2,j}=c-j(cp_3-p_2)/p_1$.
If $cp_3-p_2< 0$, $q_{2,j}$ is monotone increasing by 
Lemma \ref{keyLemma} and we see that $\cV_{\rm ns}^{(2)}(P,Q)=\emptyset$ and the
assertion follows immediately.
Assume that $cp_3-p_2\ge  0$. Then $q_{2,j}$ is  monotone decreasing for
$0\le j\le  j_1$. Thus 
$\cV_{\rm ns}^{(2)}(P,Q)\cap \cV_{\rm ns}^{(3)}(P,Q)\ne \emptyset$ if 
and only if $q_{2,j_1}=1$. If this is the case,
$Q_{j_1}$ is the unique covector in common. Thus the assertion follows from
these observations.
\end{proof}

\subsection{Normally smooth divisors  on $X_{\rm
II}$}\label{surfaceII} By using  Lemmas~\ref{numberPEi} and
\ref{numberPQ}, we can compute  the number $\rho(\Sigma ^*_{\rm
can})$.
We show this  by considering the surface  $X_{\rm II}$.
 One can 
do the same consideration for the other types of surfaces. 
Let  $X_{\rm II}: h_{\rm II}(x,y,z)=x^ay+y^b+z^c=0$.
Put  $\hat{a}:=\gcd(a,b-1), e:=\gcd(b,c)$ and
 $d:=\gcd(c(b-1),ac,ab)=e\gcd(a, c(b-1)/e)$.
The dual Newton diagram  $\Gamma ^*(h_{\rm II})_2^+$ consists of
three cones: $\Cone{P,Q}, \Cone{P,E_1}$ and $\Cone{P,E_3}$ where
 $P:=\transpose (c(b-1)/d, ac/d,ab/d)$ and
$Q:=\transpose(0,c,1)$.

The following three propositions are special cases of
 Lemmas~\ref{numberPEi} and \ref{numberPQ}.

\begin{Proposition}\label{existIPR}\sl
 $\Cone{P,E_1}$ is regular if and only if $a$ divides
$c(b-1)/e$. 
 Assume that $a\nmid (c(b-1)/e)$. Then
\begin{itemize}
\item[1)] $\cV_{\rm ns}^{(1)}(P,E_1)\ne \emptyset$ if and only if 
$ae>(b-1)c$. And in this case $\rho^{(1)}_{PE_1}
=\left[\frac{ae}{(b-1)c}\right]$.
\item[2)] $\cV_{\rm ns}^{(2)}(P,E_1)\ne \emptyset$ 
 if and only if $c|b$.
\item[3)] $\cV_{\rm ns}^{(3)}(P,E_1)\ne \emptyset$ if and only if $b|c$.
\item[4)] 
$\rho_{PE_1}=\max\{\rho_{PE_1}^{(2)},\rho_{PE_1}^{(3)},
\left[\frac{ae}{(b-1)c}\right]\}$.
  \QED
\end{itemize}
\end{Proposition}

\begin{Proposition}\label{existIPT}\sl
As $\det(P,E_3)=c\hat a/d$, $\Cone{P,E_3}$ 
is regular if and only if $d=c\hat{a}$.
 Assume that $c\hat{a}>d$.
Then
\begin{itemize}
\item[1)] $\cV_{\rm ns}^{(1)}(P,E_3)\ne \emptyset$ if and only if 
$(b-1)|a$. 
\item[2)] $\cV_{\rm ns}^{(2)}(P,E_3)\ne \emptyset$ if and only if $a|(b-1)$.
\item[3)] $\cV_{\rm ns}^{(3)}(P,E_3)\ne \emptyset$ if and only if 
$c\hat{a}>ab$ and 
$\rho^{(3)}_{PE_3}
=\left[\frac{c\hat{a}}{ab}\right]$.
\item[4)] 

$\rho_{PE_3}=\max\{\rho_{PE_3}^{(1)},\rho_{PE_3}^{(2)},
\left[\frac{c\hat{a}}{ab}\right]\}$.
\end{itemize}
\end{Proposition}
\noindent
Recall that $\rho_{P,E_i}^{(j)}\le 1$ for $i=1,3$ and $j\ne i$ 
by Lemma \ref{keyLemma}.

\begin{Proposition}\label{existIPQ}\sl
 $\Cone{P,Q}$ is regular if and only if 
$(b-1)c$ divides $ae$, or equivalently $(b-1)|a$ and $c|b\frac{a}{b-1}$.
Assume that $\Cone{P,Q}$ is not regular.
Then we have
\begin{itemize}
\item[1)] $\cV_{\rm ns}^{(1)}(P,Q)=\{Q_1\}$.
\item[2)] $\cV_{\rm ns}^{(3)}(P,Q)\ne \emptyset $ if and only if $c(b-1)>ab$.
And in this case $\rho^{(3)}_{PQ}
=\left[\frac{c(b-1)}{ab}\right]$.
\item[3)] $\cV_{\rm ns}^{(2)}(P,Q)\ne \emptyset $ if and only if
 there exist positive integers $\alpha $ and $\beta$ such that
 \begin{eqnarray}
&a\beta +d\alpha =b-1,
\quad \label{exacteq}\\
&ab\beta +d\alpha  \equiv 0\mod c(b-1).
\label{reduced-eq}\end{eqnarray}
The second condition can be replaced by 
$a\beta+1\equiv 0\quad \modulo~~ c$.
\end{itemize}
\end{Proposition}

\begin{proof}The last assertion follows from by (\ref{exacteq}) as
$ab\beta +d\alpha=(b-1)(a\beta+1)$.
\end{proof}
The  non-trivial computation is required only for $\cV_{\rm ns}^{(2)}(P,Q)$
which we  will explain more in detail.
Write $b=eb_1$ and $c=e c_1$.
\begin{Corollary}\label{non-empty} {\rm I.}
For $\cV_{\rm ns}^{(2)}(P,Q)\ne \emptyset$, it is necessary that 
\begin{eqnarray}\label{necessary-cond}
\gcd(a,c)=1,
\quad b>a,c\end{eqnarray}  
In this case, we
have
$d=e\hat a$ and 
$\cV_{\rm ns}^{(2)}(P,Q)$ is the set of covectors 
$T=(\al Q+\be P)/d$ which satisfies
\begin{eqnarray}\label{existPQ}
&a\be +e\hat a \al = b-1\label{exact-eq2}\\
& 0<\al,\be
\label{estimate-cond}\\
&b-e\hat a \al\equiv 0 ~~\modulo~c\label{cong-eq2}
\end{eqnarray} 
{\rm II.} Furthermore $\cV_{\rm ns}^{(2)}(P,Q)$ is non-empty
if  $[b/c]\ge a+\hat a$.
\end{Corollary}
\begin{proof}
From the congruence $a\be+1\equiv 0$ modulo $c$, it is clear
that 
$\gcd(a,c)=1$. Thus $d=e\gcd(a,c_1(b-1))=e\hat a$.
The equality (\ref{cong-eq2}) results from
\[a\be+1=b-d\al=e(b_1-\hat a \al)\equiv 0~~\modulo ~c\]
Thus $b>a\be\ge a$ and $b> c$. The last congruence equation
is equivalent to 
$b_1-\hat a \al\equiv 0$ $ \modulo ~c_1$.

 Assume that $[b/c]-a-\hat a\ge 0$.
As $\gcd(\hat a,b_1)=1$, there exists positive integer
$\al_0$, $0<\al_0<c_1$, such that 
$b_1-\hat a \al_0\equiv 0$ modulo $c_1$. Put $b_1-\al_0 \hat a=j_0c_1$.
We see that $j_0=b_1/c_1-\al_0\hat a/c_1> [b/c]-\hat a$.
Take $\al$ which satisfies the congruence $a \beta+1\equiv 0$
modulo $c$. Then $\alpha$ takes the form
$\al=\al_0+jc_1$ with
$j\in
\bf N$ and thus  $b_1-\hat a\al=(j_0-j\hat a)c_1$.
For the positivity of $\be$, we need to have
$0\le j<j_0/\hat a$.
The integrity of $T$ implies
\[ e(b_1-\hat a \al)-1=ec_1(j_0-j \hat a)-1
\equiv 0\quad \modulo ~ a\]
As $j$ can move $0\le j<j_0/\hat a$ and $j_0> [b/c]-\hat a\ge a$
or $j_0/\hat a>a/\hat a$,
this congruence equation has a positive solution
 $j_1,~0\le j_1\le j_0/\hat a$. 
Then put $\beta=(ec_1(j_0-j_1\hat a)-1)/a$ for such a solution $j_1$.
 This gives a
covector 
$T=(\alpha Q+\be P)\in \cV_{\rm ns}^{(2)}(P,Q)$.
\end{proof}
\begin{Example}
Consider  $X_{\rm II}: x^9y+y^b+z^8=0$ with  $b=22+36k$.
Then $e=2,\hat a=3$ and  the equation is
\begin{eqnarray*}
9\be+6\al=21+36k,\quad
9\be+1\equiv 0 ~\modulo ~8
\end{eqnarray*}
In this case, 
$[b/c]-a-\hat a=(22+36k)/8-12\ge 0$ if $k\ge 37/18$.
For $k\ge 3$ (in fact, for $k\ge 2$), we have a solution
$(\al,\be)=(6k-7,7)$.
In this case, $P=\transpose( 28+48 k, 12, 33+54k)$ and $Q=\transpose(0,8,1)$
and 
$T:=(\alpha Q+\beta Q)/(28+48k)=\transpose(7,1,8)$.
We leave the computation of the other  covectors in $\cV_{\rm ns}^{(2)}(P,Q)$ 
to the reader.
\end{Example}
\subsection{The minimality of  the  canonical toric resolutions}
\label{minimality}
We study when the canonical  toric resolution of a weighted homogeneous 
surface is minimal. Though the canonical toric resolution is not
always minimal (see Example \ref{notminimal}), we can 
expect that the minimality hold for almost all  classes of non-degenerate
surfaces. By \cite[III(6.3)]{oka}, for each weighted homogeneous  surface
 the resolution graph associated with  the canonical
toric resolution  is  star-shaped.
Hence, when the resolution graph has at least three arms, the canonical
resolution is minimal.

We have the following general statement which is very helpful 
to see if a given toric modification is minimal.

\begin{Lemma}\label{whenisminimal}\sl
Let $X:=f^{-1}(0)$ be a non-degenerate surface.
 Suppose that  $P\in \Gamma^*(f)$ is  the 
strictly  positive covector corresponding  to a  compact face $\De$
 of   the Newton boundary $\Gamma(f).$
\begin{itemize}
\item[1)] Let $\De _1, \ldots , \De_{\ell}$ be the boundary edges of $\De$.
The exceptional divisor $E(P)$ is rational if and only if
$$-\frac{6{\rm Vol}( {\rm Cone}\, \De)}{d(P;f)}+\sum_{i=1}^l(r(\De_i)+1)=2$$
where ${\rm Cone}\, \De$ is the cone over $\De$ with vertex $O$ and 
$r(\De_i)$ is the number of integral points in the interior of $\De_i$.
\item[2)]
 The canonical  toric 
resolution $\pi : \tilde{X}\longrightarrow (X,0)$ is not minimal 
if and only if there exists  a compact face $\De$  
of $\Ga (f)$ such that $E(P)$ is rational,
$E(P)^2=-1$  and $E(P)$ intersects at most two
 other exceptional divisors where $P$ is the covector corresponding
to $\De$.
\end{itemize}
\end{Lemma}
\begin{proof} The first statement is a conclusion of \cite[III(6.4)]{oka}.
The assertion 2) follows from  the 
 Castelnuovo-Enriques criterion and \cite[III \S4(A) and \S6]{oka}.
\end{proof}

\begin{Theorem}\label{isminimal}\sl
Let $X$ be one of the surfaces of type $X_{\rm II}$,
$X_{\rm III}$, $X_{\rm IV}$, 
$X_{\rm V}$, $X_{\rm VII}$ or $X_{\rm VIII}$.
We assume that $a,b,c>1$ in \ref{all-types}.
Then the canonical toric resolution of $X$ is minimal. In particular,
$\rho(X,0)=\rho(\Sigma ^*_{\rm can}).$
\end{Theorem}

\begin{proof}  
We first check
when the central exceptional divisor $E(P)$ is rational by using 
 Lemma~\ref{whenisminimal} (see also \cite[III(6.9)]{oka}).
If this is the case, we compute the number of arms from $E(P)$. If this
 number is less than 3, we show that $E(P)^2\le -2$.
Recall that the number of arms in the resolution graph is the sum of 
$r(P,Q)+1$ for non-regular cones $\Cone{P,Q}\in \Ga^*(f)_2^+$.

\vspace{.15cm}\noindent
{\rm (II)}. Let $X=X_{\rm II}: x^ay+y^b+z^c=0$.
Put $e=\gcd(b,c), \hat a=\gcd(a,b-1)$.
Then $P=\transpose(c(b-1),ac,ab)/d$ with $d=e\gcd(a,c(b-1)/e)$.
Note that $r(P,Q)+1=1$, $r(P,E_1)+1=e$ and 
$r(P,E_3)+1=\hat{a}$. By loc. cit.           
$E(P)$ is rational if and only if 1) $e=\gcd(c, a/\hat{a})=1$ or 
2) $\hat{a}=\gcd(a, c/e)=1$.
If 1) holds, then $d=\hat{a}$. We have  $\det(P,Q)=c(b-1)/\hat{a}>1$,
$\det(P,E_3)=c>1$ and  $\det(P,E_1)=a/\hat{a}$.
If $\hat{a}=a$,  $\Cone{P,E_3}$  gives $\hat a=a$ arms. 
Hence, in any case  the resolution graph of $X_{\rm II}$
has  at least three arms centered at $E(P)$.

In case 2), we have $\det(P,Q)=c(b-1)/e>1$,
 $\det(P,E_1)=a>1$ and $\det(P,E_3)=c/e$. If $e<c$, we have at least
 three arms in  the resolution graph.
Suppose that  $e=c$. Then the number of arms at $E(P)$ is $e+1\ge 3$,
unless $b=2$ and $e=c=2$. In this case,
 the resolution graph has two similar arms  and $E(P)$ is normally smooth with
$E(P)^2\le -2$.

Outline of other cases:

\vspace{.15cm}\noindent
(III) Let $X_{\rm III}: x^ay+xy^b+z^c=0$. Then 
$P=\transpose (c(b-1),c(a-1),ab-1)/d$ with $d=e\gcd(c,(ab-1)/e)$
and $e=\gcd(a-1,b-1)$. The dual Newton diagram $\Ga^*(f)_2^+$ has 3 arms $\Cone{P,E_3}$,
$\Cone{P,Q}$, $\Cone{P,R}$ where $Q=\transpose(0,c,1)$
and $R=\transpose(c,0,1)$.
The central divisor $E(P)$ is rational if and only if
$\gcd (c, (ab-1)/e)=1$. 
If  $E(P)$ is rational, then $d=e$ and 
$\det(P,Q)=c(b-1)/e>1,
\det(P,R)=c(a-1)/e>1,$ and $\det(P,E_3)=c>1$. Hence, the resolution graph
has at least three arms.

\vspace{.15cm}\noindent
(IV) Let $X_{\rm IV}: x^ay+y^bz+z^c=0$.  Then 
$P:=\transpose (bc-c+1,a(c-1), ab )/d$ with $d=e\gcd(a,(bc-c+1)/e)$
 and
$e:=\gcd (b,c-1)$. The dual Newton diagram $\Ga^*(f)_2^+$ has 3 arms $\Cone
{P,E_1}$,
$\Cone{P,Q}$, $\Cone{P,S}$ where $Q=\transpose(0,c,1)$
and $S=\transpose(1,0,a)$.
The divisor $E(P)$ is rational if and only if 
$\gcd(a,(bc-c+1)/e)=1$ which is equivalent to $d=e$. We have 
$\det (P,E_1)=a>1$, $\det (P,S)=a(c-1)/e>1$ and 
$\det(P,Q) =(bc-c+1)/e$. 
 As $\Cone{P,E_1}$ has $e$-copies of arms,  $E(P)$ has at least three arms.

\vspace{.15cm}\noindent
(V) Let $X_{\rm V}: x^ay+y^bz+z^cx=0$.
Then $P:=\transpose (bc-c+1,ca-a+1,ab-b+1)/d$ with 
 $d=\gcd (bc-c+1,ca-a+1,ab-b+1)$.
The dual Newton diagram $\Ga^*(f)_2^+$ has 3 arms
$\Cone{P,Q}$, $\Cone{P,S}$, $\Cone{P,T}$ where $Q=\transpose(0,c,1)$,
 $S=\transpose(1,0,a)$ and $T:=\transpose (b,1,0)$.
The divisor $E(P)$ is rational if and only if 
$d=1$. In this case,
we have $\det (P,Q)=bc-c+1>1$, $\det (P,S)=ca-a+1>1$ and 
$\det(P,T)=ab-b+1>1$. Thus $E(P)$ has three arms.

\vspace{.15cm}\noindent
(VII) Let $X_{\rm VII}: x^az+y^bz+z^c+t
x^{c_1}y^{c_2}=0$. Then $P=\transpose (b(c-1),a(c-1),ab)/\de$
with $\de=\gcd (b(c-1),a(c-1),ab)$.
The dual Newton diagram $\Ga^*(f)_2^+$ has 4 arms
$\Cone{P,Q}$, $\Cone{P,S}$, $\Cone{P,E_1}$, $\Cone{P,E_2}$ 
where $Q=\transpose(0,1,c_2)$ and
 $S=\transpose(1,0,c_1)$.
By the weighted homogenuity, we have the equality
$b(c-1)c_1+a(c-1)c_2= abc$
which implies that $(c-1)|ab$. Hence 
$\de=(c-1)\gcd(a,b,ab/(c-1))$.
By loc. cit., $E(P)$ is rational if and only if 
either (i) $\gcd (a, b)=\gcd (a,c-1)=1$, or 
(ii) $\gcd (a, b)=\gcd (b,c-1)=1$. By symmetry, we may assume that
 the first case (i). Then
$\de=c-1$, $\det (P,Q)=b>1$, $\det (P,S)=a>1$, $\det(P,E_1)=a>1$.
 Thus
  the resolution graph has at least three arms.

\vspace{.15cm}\noindent
(VIII)  Let $X_{\rm VIII}: x^a y+xy^b+x z^c+t
y^{c_1}z^{c_2}=0$. Then $P=\transpose (c(b-1),c(a-1),b(a-1))/\de$
with $\de=\gcd (c(b-1),c(a-1),b(a-1)).$ 
By the weighted homogenuity, we must
have
$c(a-1)c_1+b(a-1)c_2=c(ab-1)$ which implies that $
(a-1)|c(ab-1)$ and 
$cc_1+bc_2=bc+{c(b-1)}/{a-1}$. Thus 
$\de=(a-1)\gcd(b,c,c(b-1)/(a-1))$.
The dual Newton diagram $\Ga^*(f)_2^+$ has  4 arms
$\Cone{P,E_3}$, $\Cone{P,Q}$, $\Cone{P,S}$ and $\Cone{P,T}$ 
where 
$ Q=\transpose (0,c,1)$,
 $S=\transpose (c_2,0,1)$ and 
$T=\transpose (c_1,1,0)$.
The divisor  $E(P)$ is rational if and only if 
$(b-1)=k(a-1) $ for some $k\in \bf N$ and $\gcd(b,c)=1$.
Then 
$d=a-1$ and 
  $\det(P,Q)=ck>1$, $\det(P,S)=c>1$, $\det(P,T)=b>1$
and $\det(P,E_3)=c$.
 Thus the $E(P)$ has at least 3 arms.
\end{proof}

\subsection{Normally smooth divisors on $T_{p,q,r}$-surfaces}\label{pqr}
 Let $T_{p,q,r}: x^p+y^q+z^r+xyz=0$ with 
$1/p+1/q+1/r<1$. 

(1) Suppose that $p,q,r$ are pairwisely coprime and  $p<q<r$.
The diagram $\Ga ^*(f)^+_2$ has three strictly positive vertices
$P:=\transpose (rq-r-q, r,q), Q:=\transpose (r,pr-p-r,p),$ and $R:=\transpose
(q,p,pq-q-p)$. The cones $\Cone{P,E_1},\Cone{Q,E_2}$ and $\Cone{R,E_3}$ are
regular. Put $\delta:=pqr-pr-qr-pq$. Then 
$\det(P,Q)=\det(Q,R)=\det(P,R)=\delta$.

\begin{Proposition}\label{Tpqr}\sl
Under the above assumption, we have 
 $$\rho(X_{p,q,r}, O)
=\rho _{QR}^{(1)}+\rho _{QR}^{(2)}+\rho _{QR}^{(3)}
+\rho _{PR}^{(2)}+\rho _{PR}^{(3)}+\rho _{PQ}^{(3)}-2-\epsilon,
$$ where
$\varepsilon=1 $ if $p=3$, and  $\epsilon=0 $ if $p\ne 3$.
\end{Proposition}
\begin{proof} This is a summary of the following three lemmas.
\end{proof}

\begin{Lemma}\label{TpqrQR}\sl
 \begin{itemize}
\item [1)]
$\cV_{\rm ns}^{(1)}(Q,R)=\{P_k=\transpose (1,k, p-k-1)\mid p/q< k<
(rp-r-p)/r\}.$
\item [2)]
$\cV_{\rm ns}^{(2)}(Q,R)=\{P_k'=\transpose (k,1, pk-k-1)\mid r/(pr-p-r)<k<
q/p\}.$
\item [3)]
$\cV_{\rm ns}^{(3)}(Q,R)=\{P_k''=\transpose (k, pk-k-1,1)\mid  q/(pq-p-q)<k< r/p \}.$
\item[4)]  $\cV_{\rm ns}^{(1)}(Q,R)\cap \cV_{\rm ns}^{(2)}(Q,R)\cap
\cV_{\rm ns}^{(3)}(Q,R) \ne \emptyset $ if and only if $p=3$.
\item[5)] $\rho _{QR}
=\rho _{QR}^{(1)}+\rho _{QR}^{(2)}+\rho _{QR}^{(3)}-1-\epsilon$, where
$\epsilon=1 $ if $p=3$, and  $\epsilon=0 $ if $p\ne 3$.
\end{itemize}
\end{Lemma}
\begin{proof} We mainly use Theorem~\ref{existgeneral}.
Let $P':=(\be Q+\al R)/\delta
=\transpose (p_1, p_2, p_3)$.
The equation is 
$$\begin{cases}
\be r+\al q=p_1 \delta\\
\be (pr-p-r)+\al p=p_2 \delta\\
\be p+\al (pq-p-q)=p_3\delta
\end{cases}
\text{ this implies  }
\begin{cases}
\al=(pr-p-r)p_1-rp_2\\
\be =qp_2-pp_1\\
p_2+p_3=(p-1)p_1
\end{cases}$$
Hence, we have the following conclusions.

 1) $p_1=1$ if and only if there exists an  integer $p_2>0$ such that
 $\al>0$ and $\be >0$. This is equivalent to 
 $p/q<p_2<(pr-p-r)/r$. And in this case 
$P'=(1, p_2, p-1-p_2)$. 

 2) $p_2=1$ if and only if there exists an  integer $p_1>0$
such that $r/(pr-p-r) <p_1<q/p$. And in this case 
$P'=(p_1, 1, (p-1)p_1-1)$. 

3) $p_3=1$ if and only if there exists an  integer $p_1>0$ such that
$q/(pq-p-q)<p_1< r/p$. And in this case
$P'=\transpose (p_1,pp_1-p_1-1,1)$.

4) is obvious now.

5) One can see this  by comparing the three sets $\cV_{\rm ns}^{(i)}(Q,R)$.
In case $p=2$, we have  $\cV_{\rm ns}^{(1)}(Q,R)=\emptyset$ and 
$\cV_{\rm ns}^{(2)}(Q,R)\cap\cV_{\rm ns}^{(3)}(Q,R)=\{\transpose (2,1,1)\}$. Hence,
$\rho _{QR} =\rho _{QR}^{(2)}+\rho _{QR}^{(3)}-1$.

In case $p=3$, we have  $\cV_{\rm ns}^{(i)}(Q,R)\cap\cV_{\rm ns}^{(j)}(Q,R)=
\cV_{\rm ns}^{(1)}(Q,R)\cap\cV_{\rm ns}^{(2)}(Q,R)\cap\cV_{\rm ns}^{(3)}(Q,R)
=\{\transpose (1,1,1)\}$ for $i\ne j$. Hence,
$\rho _{QP} =\rho _{QR}^{(1)}+\rho _{QR}^{(2)}+\rho _{QR}^{(3)}-2$.

In case $p>3$, we have $\cV_{\rm ns}^{(1)}(Q,R)\cap\cV_{\rm ns}^{(2)}(Q,R)
=\{\transpose(1,1,p-2)\}$ and $\cV_{\rm ns}^{(1)}(Q,R)\cap\cV_{\rm ns}^{(3)}(Q,R)
=\cV_{\rm ns}^{(2)}(Q,R)\cap\cV_{\rm ns}^{(3)}(Q,R)=\emptyset$. Hence,
$\rho _{QP} =\rho _{QR}^{(1)}+\rho _{QR}^{(2)}+\rho _{QR}^{(3)}-1$.
\end{proof}

Similarly, one can prove the following two lemmas.

\begin{Lemma}\label{TpqrPR}\sl
\begin{itemize}
\item [1)]$\cV_{\rm ns}^{(1)}(P,R)=\emptyset$.
\item [2)]
$\cV_{\rm ns}^{(2)}(P,R)=\{Q'_{\ell}
=\transpose ( q -\ell-1, 1,\ell)\mid q/r< \ell < (pq-p-q)/p \}.$

\item [3)] 
$\cV_{\rm ns}^{(3)}(P,R)=\{Q''_{\ell}=\transpose (q\ell -\ell-1, \ell ,1)\mid
 p/(pq-p-q)<\ell < r/q\}.$
\item [4)] Let $Q'=\transpose (q_1, q_2, q_3)=(\be P+\al R)/\delta$.
Then $(q-1)q_2=q_1+q_3$.
\item[5)] $ \rho_{PR}=\rho _{PR}^{(2)}+\rho _{PR}^{(3)}-1.$
\QED
\end{itemize}
\end{Lemma}

\begin{Lemma}\label{TpqrPQ}\sl
\begin{itemize}
\item [1)]$\cV_{\rm ns}^{(1)}(P,Q)=\cV_{\rm ns}^{(2)}(P,Q)=\emptyset$.
\item [2)]
$\cV_{\rm ns}^{(3)}(P,Q)=\{  R'_{\ell}=\transpose (r -\ell-1, \ell ,1)\mid
r/q<\ell <(pr-p-r)/p \}$ and $\rho _{PQ}=\rho _{PQ}^{(3)}$.
\QED
\end{itemize}
\end{Lemma}

\begin{Example}
(1) Let  $p=2,q=3$ and $r\ge  7$. By the canonical subdivisions of the three
cones, we see that $\rho_{QR}=\left[\frac{r-6}{2}\right]\ge 1$,
 $\rho_{PR}=\left[\frac{r-6}{3}\right]\ge 1$, and 
$\rho_{PQ}=\left[\frac{r-3}{6}\right]$. 

\noindent
(2) Let  $p=3, q=4$ and $r>4$. By the canonical subdivisions of the three
cones, we see that $\rho_{QR}=\left[\frac{r}{3}\right]\ge 1$,
 $\rho_{PR}=\left[\frac{r}{4}\right]\ge 1$ and 
$\rho_{PQ}=[\frac{2r}{3}]-[\frac {r}{4}]-1$.
\end{Example}

(2) Another case.
 Let $f(x,y,z)=x^n+y^n+z^n+xyz$ ($n\ge 4$). 
The dual Newton diagram has three
covectors $P_i,i=1,2,3$ corresponding to the three compact faces. They are given
by 
$\transpose(n-2,1,1), \transpose(1,n-2,1),\transpose(1,1,n-2)$. And for
 $i\ne j$,
$\det (P_i,P_j)=n-3$. 
Let $B_1,\dots,B_k$ be the vertices of the canonical subdivision
of $\Cone{P_1,P_2}$ from $P_1$.
Then $B_1=(P_2+(n-4)P_1)/(n-3)=\transpose(n-3,2,1)$. Thus 
$(n-3)/(n-4)=[2,\dots,2]$  with $(n-4)$-copies of 2. 
This implies  $k=n-4$ and $B_j=\transpose(n-2-j,1+j,1),j=1,\dots, n-4$.
 In fact, by Lemma \ref{keyLemma}  the third coordinate of $B_j$ is always 1 
as both of $P_1,P_2$ have  1 as the third coordinate. Hence
$\rho_{P_1P_2}=n-4$.
The branch $\Cone{P_i,E_i}$ is regular.
Thus $\rho(V,O)=\rho(\Sigma ^*_{\rm can})=3n-9$ and
{\em  every exceptional divisor is normally  smooth.}


\section{Remarks}
\subsection{Example of the inequality 
$\rho(\Sigma^*_{\rm can})>\rho(X,O)$} Let us consider
$A_{2c-1}$-singularity,
$X=\{ x^2+y^2+z^{2c}=0\}$. The
 resolution graph has two arms and  
the central divisor $E(P)$ is a  rational 
curve with $E(P)^2=-1$. Thus  we have
to blow-down the central divisor once ( Example (6.7.1) in
\cite[III]{oka} ).  However in this example,
the central exceptional divisor is not normally smooth, i.e., 
the extra blowing-up is line-admissible. So 
 $\rho(\Sigma^*_{\rm can}) =\rho(X,O)$.
The following gives an example of  $\rho(\Sigma^*_{\rm can}) >\rho(X,O)$.
\begin{Example}\label{notminimal}\rm
Let $X$ be defined by  $h=xy+y^{bc}+z^c$ with $b,c\ge 2$. 
This is an $A_{c-1}$-singularity
and  a special case of
 $X_{\rm II}$ with 
$ P:=\transpose (bc-1,1,b)$
and $ Q:=\transpose(0,c,1).$

Since $\det(P,E_1)=\det(P,E_3)=1$ and $\det(P,Q)=bc-1$, 
we make the canonical subdivision of
 $\Cone{P,Q}$. 
The first covector $T_1$ from $P$ is given by 
\[T_1=(Q+(bc-c-1)P)/(bc-1)=\transpose
(bc-c-1,1,b-1)\]
We have the continuous fraction expansion
$(bc-1)/(bc-c-1)=[2,\dots,2,3,2,\dots, 2]$
where the number of $2$ in  the first 2-series
(respectively in the second $2$-series)
 is $(b-2)$ (resp. $c-2$).
Thus we have $c+b-3$ covectors
$T_1,\dots, T_{b+c-3}$. The
exceptional divisor  $E(P)$ is  rational with
$E(P)^2=-1$  and $E(T_j)$
with self intersection number
$E(T_j)^2=-2 $ for $j\ne b-1$ and $-3$ for $j=b-1$
(see Theorem (6.3), Chapter III, \cite{oka}).
In fact first $b-2$ covectors are given by
\begin{eqnarray*}
&Q_j=\transpose (cb-jc-1,1,b-j),\quad j=1,\dots,b-1\\
& Q_{b-1+j}=\transpose
(c-j-1,j+1,1),\quad j=1,\dots, c-2
\end{eqnarray*}
and we see that they are normally minimal.
To get a minimal reslution, we need to blow down $b-1$ divisors
$E(P),E(T_1),\dots, E(T_{b-2})$ in this order. Then the 
self-intersection number of $E(T_{b-1})$ changes to 
 $-2$ and we
get 
$A_{c-1}$ graph. In this example,
we have $\rho(X,O)=c-1$
and $\rho(\Sigma_{\rm can}^*)=b+c-2$.
\end{Example}

\subsection{Parametrization of lines}
The normally smooth divisors
on a surface $X$ correspond to the lines on $X$. By using a toric resolution,
 one can give the exact parameterizations of the lines on $X$.
This was done already for the Pham-Brieskorn surfaces in
\cite{JOPS}. 
\begin{Proposition}
Suppose that we have a line $L$ in a
non-degenerate surface $X: f(x,y,z)=0$ and  assume that $L$ is
parametrized as
\[x(t)=\al t^a+\al_1t^{a+1}\dots,
\quad  y(t)=\be t^b+\be_1 t^{b+1}+\dots, 
\quad z(t)=\ga t^c+\ga_1 t^{c+1}+\dots\] with $\al, \be, \ga\ne
0$ and
$\min(a,b,c)=1$. Let $P=\transpose (a,b,c)$. Then the  pull
back of $L$ intersects $E(P)$ transversally 
and $f_P(\al,\be,\ga)=0$. 
Conversely any curve in ${\cal L}_{E(P)}$
has such a parametrization.
\end{Proposition}
\begin{Example} \label{EII-abb}
(1)  Let $X$ be defined by  $h=x^ay+y^b-z^b=0$ with 
$a=a_1(b-1)$ and $a_1>1$.
This is a special case of $X_{\rm II}$. We use the notations in
 \S\ref{surfaceII}.
Note that
$P=\transpose (1, a_1, a_1 )$, $Q=\transpose(0,b,1)$,
$\det(P,Q)=\det(P,E_3)=1$  and $\det(P,E_1)=a_1$. By canonical
 subdivision of
$\Cone{P,E_1}$ we have $R_i:=\transpose (1,i,i)$ with
 $i=0,1, \ldots , i_1=a_1$,
where $R_0:=E_1$ and $R_{i_1}:=P$. Hence $\rho _{PE_1}=a_1-1$.
 Since $r(P,E_1)+1=b$, 
each $E(R_i)$ has  $b$ components.
By \cite[III(6.3)]{oka}, $E(P)^2=-b<-1$. Hence $\pi$ is minimal and 
$\rho (X,0)=b(a_1-1)+1$.
The restriction of $\pi$ on the toric chart associated with 
 $\sigma _i:={\rm Cone}~(R_i,R_{i-1},E_2)$ is  given by
$$\pi _{\sigma _i}:\quad x=uv, \quad y=u^iv^{i-1}w, \quad z=u^iv^{i-1}.$$
and the pull-back of $h$ is given by  
$$h\circ \pi _{\sigma _i}=u^{ib}v^{(i-1)b}
\left(u^{(a_1-i)(b-1)}v^{(a_1-i+1)(b-1)}w+w^b-1\right)$$
The divisor
$E(R_i)$ is defined by $u=0$ and $w^b-1=0$, hence $E(R_i)$ has $b$ 
components.
 On this toric  chart, the resolution $\tilde X$ of $X$ is
defined by 
\[\tilde h_i(u,v,w):=u^{(a_1-i)(b-1)}v^{(a_1-i+1)(b-1)}w+w^b-1=0\]
 and in a
neighborhood
 of $q\in E(R_i)$ we take $u,v$ to 
 be the local coordinates of $\tilde X$. Let $q=(0,s)$ in this
coordinates.
We consider the lines $ C_s$ 
defined by $t\mapsto (t,s)$.
 The image of $C_s$ by $\pi _{\sigma _i}$ is given by
$$\pi_{\sigma _i}(C_s):\quad x=st,\quad y=s^{i-1}w_k(t,s)t, \quad
z=s^{i-1}t^i,$$
where $w_k(t,s)$  is  the  solution of $\tilde h_i(t,s,w)=0$
with $w_k(0)=\exp(2\pi ki/b)$.
As a special case, take $i=1$. Then $C_s$ is a normal line
on $E(Q_1)$. When we moves $s\to 0$, this line approaches to 
$E(E_1)$ and  $w_k(t)\equiv \exp (2 k\pi i/b)$
and   the image  is the obvious line $t\to (x,y,z)=(0,w_k t,t)$.

(2) Let $X=T_{2,3,7}: x^2+y^3+z^7+xyz=0$. We have three covectors
\[P=\transpose(11,7,3),~ Q=\transpose(7,5,2),~R=\transpose (3,2,1)\]
and we do not need any other covector. Consider the toric chart
$\si:=(Q,R,E_3)$ with coordinates $(u,v,w)$.
Then the line $u=1,v=t$ produces a line
parametrized as
$t\mapsto (t^3,t^2, -2t+128 t^2+\dots)$.
\end{Example}

\subsection{Obvious lines on surfaces}
We consider a surface $X=\{f(x,y,z)=0\}$
where $f$ has a non-degenerate Newton boundary.
There are  surfaces having  obvious lines which can be read off from the
polynomial defining the surface.

\noindent
(1) Assume that $f(x,y,z)$ is not convenient and assume for example
$\{y=z=0\}\subset X$. Then as we have seen in Lemma~\ref{non-compact-face},
there is a unique non-compact face, different from the  coordinate planes,
which has the covector of the type
$Q=\transpose (0,c,1)$ or $\transpose (0,1,c)$ and a unique covector
$P$ such that $\Cone{P,Q}$ is in $\Ga^*(f)_2^{+}$ and $P$
corresponds to a compact face. Let $Q_1,\dots, Q_k$ be the
covectors defining the  canonical regular subdivision from $Q$. Then $Q_1$ 
is a normally smooth divisor and $\cal L_{\rm Q_1}$ contains the canonical line
$\{y=z=0\}$.

\noindent
(2)  Assume that  $h(x,y):=f(x,y,0)$
(the section of $f$ with $z=0$) is a
non-monomial  homogeneous
polynomial of degree
$d$. Then we can factor $h(x,y)=c x^a y^b \prod_{i=1}^k (y-\al_i x)$.
Thus $X$ has the lines $z=0,y=\al_i x$ for $i=1,\dots,k$.
Combinatorially this says the following.
 There exists a compact  face $\De$ such that 
$\De\supset \De(h)$. The corresponding covector takes the form
$P=\transpose (p,p,r)$ with 
$\gcd(p,r)=1$. Then the first covector
$Q_1$ from $E_3$ in the canonical regular subdivision of $\Cone{P,E_3}$ 
takes the form $Q_1=\transpose ( 1,1,s)$ 
with $s=1+ [r/p]$. So we can see that $Q_1\in \cV_{\rm ns}(P,E_3)$.
A typical example is $T_{n,n,n}: x^n+y^n+z^n-xyz=0$. 
Another example is (1) of Example \ref{EII-abb}.

\noindent
(3)  Assume that the monomial $x^A$ in $f$ such that  $(A,0,0)\in\Ga(f)$.
We say that $x^A$ is negligibly truncatable if 
$f_t(x,y,z)=(f(x,y,z)-f(x,0,0))+ tf(x,0,0)$ defines a $\mu$-constant
family for $0\le t\le 1$ (cf. \cite{oka2}). Assume for example,
the monomials $x^ay$ and $x^bz^c$ are on the non-compact face of $\Ga(f_0)$.
Let $Q':=\transpose(c/d,c(A-a)/d,(A-b)/d)$ with $d=\gcd(c,A-b)$. 
The covector $Q'$ corresponds to 
the negligible compact face of $f_1$ containing 
$(a,1,0),(b,0,c)$ and $(A,0,0)$.
Then there is a normally smooth
divisor on  $\Cone{Q,E_3}$. In fact, $\det(Q',E_3)=c/d$. If 
$c=d$, $Q$ gives normally smooth divisor. If $c>d$, the first covector
$Q_1'$ of the canonical regular subdivision of $\Cone {Q',E_3}$
is normally smooth.
An example is given by $f(x,y,z)=x^2y+y^2+z^5+x^5$. Then $x^5$ is negligibly truncatable.

\noindent
(4) Assume that $\Ga(f)$ has a compact face whose  covector $P$
 has 1 in its coefficients. Then $E(P)$ is a normally smooth divisor.
This is the case, for example, if $P=\transpose (1,1,1)$ and
$f_P(x,y,z)$ has a
two-dimensional  support. We can see easily that
$E(P)$ is isomorphic to the projective curve $f_P(x,y,z)=0$ in ${\mathbb P}^2$.
The tangent cone of $X$ at $O$ is given by the cone of $f_P=0$.

\subsection{Normally smooth divisors on complete intersections}
In this paper we mainly considered  normally smooth divisors on two
dimensional hypersurface singularities. However
every assertion   can be generalized to non-degenerate complete
intersections. We give an example.
Consider the surface given by 
$X=\{f_1(x,y,z,w)=f_2(x,y,z,w)=0\}$ where $f_1$ and $f_2$ 
has the same Newton boundary.
 Assume that $f_1,f_2$ are Pham-Brieskorn  
polynomials of the same type, with generic coefficients:
\[f_i=a_i x^{p_1}+b_i y^{p_2}+c_i z^{p_3}+d_i w^{p_4},i=1,2\]
We assume that $p_1,\dots,p_4\ge 2$ and mutually coprime.
Then  the dual Newton diagram $\Ga^*(f_1,f_2)$
is the same with  $\Ga^*(f_i)$ and $\Ga^*(f)_2^+$  is star-shaped with 
 the center \newline
$P=\transpose(p_2p_3p_4,p_1p_3p_4,p_1p_2p_4,,p_1p_2p_3)$ and 
four arms $\Cone{P,E_i},i=1,\dots,4$.
We consider the $\Cone{P,E_1}$. First $\det(P,E_1)=p_1$.
By Lemma \ref{existspecial}, $\cV_{ns}^{(i)}(P,E_1)=\emptyset$ 
for $2\le i\le 4.$
As for $\cV_{ns}^{(1)}(P,E_1)\ne \emptyset$ if and only if 
$p_2p_3p_4<p_1$ and putting $r=[p_1/p_2p_3p_4]$,
$\cV_{ns}^{(1)}(P,E_1)=\{Q_j=(j P+ (p_1-j p_2p_3p_4)E_1)/p_1;j=1,\dots,r\}$.


\end{document}